\documentclass[onecolumn,reprint,a4paper,superscriptaddress,prb]{revtex4}

\usepackage{amsmath}    
\usepackage{verbatim}   
\usepackage{hyperref}  
\usepackage{bbold}
\usepackage{mathrsfs}
\usepackage{textcomp}
\usepackage{amsthm}
\usepackage{amssymb}
\usepackage{scrextend}
\usepackage{titlesec}
\allowdisplaybreaks

\titleformat*{\subsubsection}{\large\bfseries}

\usepackage{scrextend}
\changefontsizes[15pt]{11pt}

 \usepackage[
    top    = 2.00cm,
    bottom = 2.00cm,
    left   = 1.75cm,
    right  = 1.75cm]{geometry}

\begin{document}

\title{Two-dimensional Fourier transformations and double Mordell integrals} 

\author{Martin Nicholson}

\begin{abstract}  
      Fourier transformations of several functions of one and two variables are evaluated and then used to derive some integral and series identities. It is shown that certain double Mordell integrals can be reduced to a sum of products of one-dimensional Mordell integrals. As a consequence of this reduction, a quadratic polynomial identity is found connecting products of certain one-dimensional Mordell integrals. An integral that depends on one real valued parameter is calculated reminiscent of an integral previously calculated by Ramanujan and Glasser. Some connections to elliptic functions and lattice sums are discussed.
\end{abstract}

\maketitle

\section{Introduction: self-reciprocal Fourier transformations}

Define the cosine and sine Fourier transformations
by the usual formulas
\begin{equation}\label{cos}
    f_c(t)=\sqrt{\frac{2}{\pi}}\int_0^\infty f(x)\cos tx\,dx,
\end{equation}
\begin{equation}\label{sin}
    f_s(t)=\sqrt{\frac{2}{\pi}}\int_0^\infty f(x)\sin tx\,dx.
\end{equation}
Functions that are equal to their own cosine Fourier transform, i.e. that satisfy the equation $f(x)=f_c(x)$, are called self-reciprocal functions of the first kind, and functions that are equal to their own sine Fourier transform $f(x)=f_s(x)$, are called self-reciprocal functions of the second kind[\onlinecite{hardy}]. Some examples of the functions of the first kind include
\begin{equation}\label{cosineselfreciprocal}
\frac{1}{\cosh\sqrt{\frac{\pi}{2}}x},~ \frac{\cosh \frac{\sqrt{\pi}x}{2}}{\cosh \sqrt{\pi}x},~\frac{1}{1+2\cosh \sqrt{\frac{2\pi}{3}}x},~\frac{\cosh\frac{\sqrt{3\pi}x}{2}}{2\cosh \sqrt{\frac{4\pi}{3}} x-1},~\frac{\cosh\sqrt{\frac{3\pi}{2}}x}{\cosh \sqrt{2\pi}x-\cos\sqrt{3}\pi}.
\end{equation}
And here are some functions of the second kind
\begin{equation}\label{sineselfreciprocal}
    \frac{\sinh \frac{\sqrt{\pi}x}{2}}{\cosh \sqrt{\pi}x}, \frac{\sinh\sqrt{\frac{\pi}{6}}x}{2\cosh \sqrt{\frac{2\pi}{3}}x-1}, \frac{\sinh \sqrt{\frac{2\pi}{3}}x}{\cosh \sqrt{\frac{3\pi}{2}}x}, \frac{\sinh\sqrt{\pi}x}{\cosh \sqrt{2\pi} x-\cos\sqrt{2}\pi}.
\end{equation}
 The first three functions of \eqref{cosineselfreciprocal} and the first two functions of \eqref{sineselfreciprocal} were known to Ramanujan and their detailed study can be found in the book [\onlinecite{berndt}]. The third function in \eqref{sineselfreciprocal} is taken from the article [\onlinecite{cais}] where many other hyperbolic self reciprocal functions are given along with a general method for generating them. The last two functions in \eqref{cosineselfreciprocal} and the last function in \eqref{sineselfreciprocal} appear to be new. One can show that \eqref{cosineselfreciprocal} are the only self reciprocal functions of the form $\frac{\cosh \alpha x}{\cosh x+c}$.

There is a well known general recipe to find self reciprocal functions ([\onlinecite{titchmarsh}, ch. 9]). Since $(f_c)_c=f$, the sum
$$
f(x)+f_c(x)
$$
is a self-reciprocal function of the first kind for an arbitrary function $f(x)$. Obviously this approach works also for functions of the second kind.

It might seem that this settles the question of finding all self-reciprocal functions completely. However this is not so because this approach is not helpful in finding interesting particular self-reciprocal functions. It is much more gratifying to now that the functions in \eqref{cosineselfreciprocal} are self-reciprocal as opposed to knowing that the function
$$
e^{-x}+\sqrt{\frac{2}{\pi}}\frac1{1+x^2}
$$
is self-reciprocal. A more useful general theory suitable for these purposes of finding particular transformations has been developed by Goodspeed, Hardy and Titchmarsh (see [\onlinecite{titchmarsh}] for a nice account of this theory).

One might ask, what are these particular transformations useful for? The answer is they lead to some interesting integral and series transformation formulas, among other things. For example, Hardy and Ramanujan [\onlinecite{hardy},\onlinecite{ramanujan}] used self reciprocal functions to obtain transformation formulas such as
\begin{equation}\label{int2}
\sqrt{\alpha}\int\limits_0^\infty\frac{\cosh \frac{\alpha x}{2}}{\cosh \alpha x}\,e^{-x^2}dx= \sqrt{\beta}\int\limits_0^\infty\frac{\cosh \frac{\beta y}{2}}{\cosh \beta y}\,e^{-y^2}dy,\qquad \alpha\beta=2\pi,
\end{equation}
\begin{equation}\label{int3}
\sqrt{\alpha}\int\limits_0^\infty\frac{\sinh \frac{\alpha x}{2}}{\sinh \alpha x}\,xe^{-x^2}dx= \sqrt{\beta}\int\limits_0^\infty\frac{\sinh \frac{\beta y}{2}}{\sinh \beta y}\,ye^{-y^2}dy,\qquad \alpha\beta=2\pi.
\end{equation}

Another type of identities are obtained by application of the Poisson summation formula, which for an even function $\phi(x)$ can be stated in the symmetric form [\onlinecite{titchmarsh}]
\begin{equation}\label{poisson1}
    \sqrt{\alpha}\sum_{n=-\infty}^\infty \phi(\alpha n)=\sqrt{\beta}\sum_{n=-\infty}^\infty \phi_c(\beta n),\qquad \alpha\beta=2\pi.
\end{equation}
Similarly, for an odd function $\psi(x)$
\begin{equation}
    \sqrt{\alpha}\sum_{n=1}^\infty\chi(n) \psi(\alpha n)=\sqrt{\beta}\sum_{n=1}^\infty\chi(n) \psi_s(\beta n),\qquad \alpha\beta=\frac{\pi}{2},
\end{equation}
where $\chi(n)=\sin\frac{\pi n}{2}$ is a primitive character of modulus 4. For example, application of \eqref{poisson1} to the first function in \eqref{cosineselfreciprocal} gives
\begin{equation}\label{elliptic}
    \sqrt{\alpha}\sum_{n=-\infty}^\infty\frac{1}{\cosh \pi \alpha n}=\sqrt{\beta}\sum_{n=-\infty}^\infty\frac{1}{\cosh \pi \beta n}, \qquad \alpha\beta=1.
\end{equation}
Let $q=e^{-\pi\alpha}$ be the base of elliptic functions with modulus $k$,$~k'=\sqrt{1-k^2}$ the complementary modulus and $K=K(k)$,$~K'=K(k')$ the complete elliptic integrals of the first kind. Then [\onlinecite{ww}, ch. $22.6$] $q'=e^{-\pi\beta}$ is the base of elliptic functions with modulus $k'$ and
\begin{equation}\label{K}
    K=\frac{\pi}{2}\sum_{n=-\infty}^{\infty}\frac{1}{\cosh \pi\alpha n}.
\end{equation}
So \eqref{elliptic} is nothing but $q=e^{-\pi\frac{K'}{K}}$ in the more familiar notation of the theory of elliptic functions.

Functions \eqref{cosineselfreciprocal},\eqref{sineselfreciprocal} imply certain symmetric relations for the Lerch zeta function ([\onlinecite{berndt}], ch. $18.5$). For example the fourth function in \eqref{sineselfreciprocal} leads to the identity
\begin{equation}\label{lerch}
    \sum_{n=-\infty}^\infty\frac{\sin\frac{\sqrt{2}\pi n}{p}}{{\left| n+\frac{p}{\sqrt{2}}\right|}^{\frac12}}=
\sum_{n=-\infty}^\infty\frac{\sin\frac{\sqrt{2}\pi n}{q}}{{\left| n+\frac{q}{\sqrt{2}}\right|}^{\frac12}},\qquad pq=1,\ \frac{1}{\sqrt{2}} < p < \sqrt{2}.
\end{equation}

\section{Functions of two variables}

One may also consider self reciprocal Fourier functions of two variables. Apart from the non-interesting factorizable functions of this form there are quite non-trivial functions. To find some of them we use the following observation:
If $f(x,y)=f(y,x)$ and $$\sqrt{\frac{2}{\pi}}\int_0^\infty f(x,y)\cos ax\, dx=g(a,y)=g(y,a),$$ (in other words, if partial Fourier transform of a symmetric function is symmetric) then $f(x,y)$ is a self-reciprocal Fourier function of two variables, i.e. 
$$
{\frac{2}{\pi}}\int\limits_0^\infty \int\limits_0^\infty f(x,y)\cos ax\cos by \,dxdy=f(a,b).
$$
$\it{Example:}$ Since ([\onlinecite{GR}], formula 3.981.8)
$$
\int\limits_0^\infty\frac{\sin xy}{\sinh\sqrt{\pi} x}\cos ax
\,dx=\frac{\sqrt{\pi}}{2}\frac{ \sinh\sqrt{\pi}y}{\cosh\sqrt{\pi}y+\cosh\sqrt{\pi}a}$$
we get a pair of self-reciprocal Fourier transformations 
\begin{equation}\label{cos1}
    \frac{2}{\pi}\int\limits_0^\infty \int\limits_0^\infty\frac{\cos ax\cos by}{\cosh\sqrt{\pi}x+\cosh\sqrt{\pi}y}\,dxdy=\frac{1}{\cosh\sqrt{\pi}a+\cosh\sqrt{\pi}a},
\end{equation}
\begin{equation}\label{cos2}
    \frac{2}{\pi}\int\limits_0^\infty \int\limits_0^\infty\frac{\sin xy}{\sinh\sqrt{\pi} x\,\sinh\sqrt{\pi}y}\cos ax\cos by\,dxdy=\frac{\sin ab}{\sinh\sqrt{\pi} a\, \sinh\sqrt{\pi}b}.
\end{equation}
Though not a self reciprocal function, note the curious transformation
\begin{equation}\label{f1}
    {\frac{2}{\pi}}\int\limits_0^\infty \int\limits_0^\infty \frac{\cos xy}{\cosh\sqrt{\frac{\pi}{2}} x\, \cosh\sqrt{\frac{\pi}{2}} y}\cos ax\cos by\,dxdy=\frac{\sin ab}{\sinh\sqrt{\frac{\pi}{2}} a \,\sinh\sqrt{\frac{\pi}{2}} b}.
\end{equation}
More self-reciprocal functions of one and two variables can be found in [\onlinecite{blog}].

Poisson summation formula \eqref{poisson1} is easily generalized to even functions of two variables as follows
\begin{equation}\label{poisson2}
    \sqrt{\alpha\beta}\sum_{m,n=-\infty}^\infty \phi(\alpha m,\beta n)=\sqrt{\gamma\delta}\sum_{m,n=-\infty}^\infty \phi_c(\gamma m,\delta n),\qquad \alpha\gamma=\beta\delta=2\pi,
\end{equation}
where
\begin{equation*}
    \phi_c(t,s)=\frac{2}{\pi}\int\limits_0^\infty\int\limits_0^\infty \phi(x,y)\cos tx\cos sy\,dxdy.
\end{equation*}

It is instructive to see what happens if \eqref{poisson2} is applied to \eqref{f1}. Straightforward calculation shows that
\begin{equation*}
    \sqrt{\alpha\beta}\sum_{m,n=-\infty}^\infty \frac{\cos \alpha\beta m n}{\cosh\sqrt{\frac{\pi}{2}} \alpha m\cdot \cosh\sqrt{\frac{\pi}{2}} \beta n}=\sqrt{\gamma\delta}\sum_{m,n=-\infty}^\infty \frac{\sin \gamma\delta m n}{\sinh\sqrt{\frac{\pi}{2}} \gamma m \cdot\sinh\sqrt{\frac{\pi}{2}} \delta n},\qquad \alpha\gamma=\beta\delta=2\pi.
\end{equation*}
Here it is assumed that the terms with $m=0$ or $n=0$ on the RHS of are understood as the limits $\lim_{m\to 0}$,$~\lim_{n\to 0}$. Setting $\delta=\alpha,\ \gamma=\beta$ and making the replacement $\alpha\to \sqrt{2\pi}\alpha$, $\beta\to \sqrt{2\pi}\beta$ one obtains (care should be taken to simplify the sum on the right)
\begin{equation}\label{Legendre}
    \sum_{n=-\infty}^\infty\frac{1}{\cosh \pi n\alpha}\sum_{n=-\infty}^\infty\frac{1}{\cosh \pi n\beta}=\frac{2}{\pi}+4\sum_{n=1}^\infty\frac{\alpha n}{\sinh \pi n\alpha}+4\sum_{n=1}^\infty\frac{\beta n}{\sinh \pi n\beta},\quad \alpha\beta=1.
\end{equation}
It is known that [\onlinecite{ww}, ch. $22.735$]
$$
\sum_{n=1}^\infty\frac{n}{\sinh \pi n\alpha}=\frac{K(K-E)}{\pi^2},
$$
with the same notations as in \eqref{K} and $E=E(k)$ complete elliptic integral of the second kind. Therefore \eqref{Legendre} is Legendre's relation $EK'+E'K-KK'=\frac{\pi}{2}$ in disguise.

Hyperbolic functions provide many other transformations. Let's start with the calculation of the integral
$$
J=\int\limits_0^\infty \int\limits_0^\infty\frac{\cos xy}{\cosh p x\cosh\frac{\pi y}{p}}\cos ax\cos by\,dxdy.
$$
By formula $3.981.10$ from [\onlinecite{GR}]: 
\begin{align*}    
   J&=\int\limits_0^\infty\frac{\cos by}{\cosh\frac{\pi y}{p}}\, dy\int\limits_0^\infty\frac{\cos xy}{\cosh p x}\cos ax\,dx\\
&=\int\limits_0^\infty\frac{\cos by}{\cosh\frac{\pi y}{p}} \cdot\frac{\pi}{p}\frac{\cosh\frac{\pi a}{2p}\cosh\frac{\pi y}{2p}}{\cosh\frac{\pi a}{p}+\cosh\frac{\pi y}{p}}\,dy\\
&=\frac{\pi}{p}\frac{\cosh\frac{\pi a}{2p}}{\cosh\frac{\pi a}{p}}\cdot\int\limits_0^\infty\left(\frac{\cosh\frac{\pi y}{2p}}{\cosh\frac{\pi y}{p}} -\frac{\cosh\frac{\pi y}{2p}}{\cosh\frac{\pi a}{p}+\cosh\frac{\pi y}{p}}\right)\cos by\,dy\\
&=\frac{\pi}{\sqrt{2}}\cdot \frac{\cosh{\frac{\pi a}{2p}\cosh \frac{p b}{2}}}{\cosh\frac{\pi a}{p}\cosh pb}-\frac{\pi}{2}\cdot\frac{\cos ab}{\cosh \frac{\pi a}{p}\cosh pb},
\end{align*}
so finally
\begin{equation}\label{I1}
    \frac{2}{\pi}\int\limits_0^\infty \int\limits_0^\infty\frac{\cos xy}{\cosh p x\cosh\frac{\pi y}{p}}\cos ax\cos by\,dxdy=\sqrt{2}\cdot \frac{\cosh{\frac{\pi a}{2p}\cosh \frac{p b}{2}}}{\cosh\frac{\pi a}{p}\cosh pb}-\frac{\cos ab}{\cosh \frac{\pi a}{p}\cosh pb}.
\end{equation}
We see that the right hand side is the original function (taken with the minus sign) up to an additional term, which a factorizable function. 

Applying Poisson summation \eqref{poisson2} to \eqref{I1} one finds 
\begin{equation}\label{landen}
    \sqrt{2}\sum_{m=-\infty}^\infty\frac{1}{\cosh\pi\alpha m}\sum_{n=-\infty}^\infty\frac{1}{\cosh\pi\beta n}=\sum_{m=-\infty}^\infty\frac{\cosh\frac{\pi\alpha m}{2}}{\cosh\pi\alpha m}\sum_{n=-\infty}^\infty\frac{\cosh\frac{\pi\beta n}{2}}{\cosh\pi\beta n},\qquad \alpha\beta=2.
\end{equation}
\eqref{landen} is equivalent to the modulus transformation of Landen's transform, i.e. $(1+k_1)(1+k')=2$ in the notation of the book [\onlinecite{ww}]. Indeed, if $\alpha=\frac{K(k')}{K(k)}$, $~\beta=\frac{\Lambda(k_1)}{\Lambda(k_1')}$, then
$$
\Lambda'=\frac{\pi}{2}\sum_{n=-\infty}^\infty\frac{1}{\cosh\pi\beta n},
$$
$$
\text{dn}\left(\tfrac{iK'}{2},k\right)=\frac{\pi}{2K}\sum_{n=-\infty}^\infty\frac{\cosh\frac{\pi\alpha n}{2}}{\cosh\pi\alpha n},
$$
$$
\text{dn}\left(\tfrac{i\Lambda}{2},k_1'\right)=\frac{\pi}{2\Lambda'}\sum_{n=-\infty}^\infty\frac{\cosh\frac{\pi\beta n}{2}}{\cosh\pi\beta n}.
$$
Since $\text{dn}\left(\frac{iK'}{2},k\right)=\sqrt{1+k}$, eq. \eqref{landen} reduces to $(1+k)(1+k_1')=2$, as required.

There is an integral analogous to \eqref{I1} involving odd functions:
\begin{equation}\label{I2}
    \frac{2}{\pi}\int\limits_0^\infty \int\limits_0^\infty\frac{\sin xy}{\cosh p x\cosh\frac{\pi y}{p}}\sin ax\sin by\,dxdy=\sqrt{2}\cdot \frac{\sinh{\frac{\pi a}{2p}\sinh \frac{p b}{2}}}{\cosh\frac{\pi a}{p}\cosh pb}-\frac{\sin ab}{\cosh \frac{\pi a}{p}\cosh pb}.
\end{equation}
Just to illustrate what kind of transformations one can get by considering more complicated functions:
\begin{align*}
    \frac{4}{\pi}\int\limits_0^\infty \int\limits_0^\infty\frac{\cos xy \cos ax\cos by\,dxdy}{(1+2\cosh x)(1+2\cosh \frac{2\pi y}{3})}=\sqrt{3}\sin ab\,\frac{\cosh\frac{b}{2}}{\sinh\frac{3b}{2}}\frac{\cosh\frac{\pi a}{3}}{\sinh\pi a}-\frac{1+\cos ab}{(1+2\cosh \frac{2\pi a}{3})(1+2\cosh b)},
\end{align*}
\begin{align*}
    \frac{4}{\pi}\int\limits_0^\infty \int\limits_0^\infty\frac{\sin xy \sin ax\sin by\,dxdy}{(1+2\cosh x)(1+2\cosh \frac{2\pi y}{3})}=\frac{\sqrt{3}(1-\cos ab)\cosh\frac{b}{2}\cosh\frac{\pi a}{3}}{\sinh\frac{3b}{2}\sinh\pi a}-\frac{\sin ab}{(1+2\cosh \frac{2\pi a}{3})(1+2\cosh b)}.
\end{align*}

\section{Case studies of several two-dimensional Mordell integrals}

Let's multiply \eqref{I1} by $e^{-({a^2}+{b^2})/2}$ and integrate with respect to $a$ and $b$
\begin{align*}
    &\int\limits_0^\infty \int\limits_0^\infty\frac{\cos xy}{\cosh p x\cosh\frac{\pi y}{p}}\,e^{-({x^2}+{y^2})/2}\,dxdy=\\
&\sqrt{2}\cdot 
\int\limits_0^\infty \int\limits_0^\infty\frac{\cosh{\frac{\pi a}{2p}\cosh \frac{p b}{2}}}{\cosh\frac{\pi a}{p}\cosh pb}\,e^{-({a^2}+{b^2})/2}\,dadb-
\int\limits_0^\infty \int\limits_0^\infty\frac{\cos ab}{\cosh \frac{\pi a}{p}\cosh pb}\,e^{-({a^2}+{b^2})/2}\,dadb.
\end{align*}
This can be written in the following symmetrical form
\begin{equation}\label{landen2}
    \sqrt{2}\cdot\int\limits_0^\infty \int\limits_0^\infty\frac{\cos 2xy}{\cosh \alpha x\cosh\beta y}\,e^{-x^2-y^2}dxdy=\int\limits_0^\infty\frac{\cosh \frac{\alpha x}{2}}{\cosh \alpha x}\,e^{-x^2}dx\cdot \int\limits_0^\infty\frac{\cosh \frac{\beta y}{2}}{\cosh \beta y}\,e^{-y^2}dy,\qquad \alpha\beta=2\pi.
\end{equation}
Note the similarity of \eqref{landen2} with the Landen transform \eqref{landen}. Since Mordell integrals can be understood as continous analogs of theta functions [\onlinecite{berndt}], \eqref{landen2} can be understood as Landen's transform for Mordell integrals. However the factorization on the left side of \eqref{landen2} does not occur because of the function $\cos 2xy$ in the integrand (in the discrete case it was possible to choose the parameters so that $\cos 2xy$ didn't have any mixing effect on the two series, so the double series factorized; unfortunately this is not possible for an integral).

Combining \eqref{landen2} with \eqref{int2} leads to
\begin{equation}\label{factorization1}
    \int\limits_0^\infty \int\limits_0^\infty\frac{e^{-x^2-y^2}\cos 2xy}{\cosh \alpha x\cosh({2\pi y}/{\alpha})}\,dxdy=\frac{\alpha}{2\sqrt{\pi}}\left(\int\limits_0^\infty\frac{\cosh \frac{\alpha x}{2}}{\cosh \alpha x}\,e^{-x^2}dx \right)^2.
\end{equation}
${\it{Corollary}}~1.$ 
$$
\int\limits_0^\infty \int\limits_0^\infty\frac{\cos\frac{\pi}{2}
\left(nx^2-\frac{y^2}{n}\right)\cos \pi xy}{\cosh \pi x\cosh \pi y}\,dxdy=\frac{\sqrt{n}}{2}I_1^2-\frac{\sqrt{n}}{2}I_2^2+\sqrt{n}I_1I_2,
$$
$$
\int\limits_0^\infty \int\limits_0^\infty\frac{\sin\frac{\pi}{2}
\left(nx^2-\frac{y^2}{n}\right)\cos \pi xy}{\cosh \pi x\cosh \pi y}\,dxdy=\frac{\sqrt{n}}{2}I_2^2-\frac{\sqrt{n}}{2}I_1^2+\sqrt{n}I_1I_2,
$$
where $\ \displaystyle I_1=\int\limits_0^\infty\frac{\cosh \frac{\pi x}{2}}{\cosh \pi x}\cos{\frac{\pi n x^2}{2}}\,dx$,$\displaystyle\ I_2=\int\limits_0^\infty\frac{\cosh \frac{\pi x}{2}}{\cosh \pi x}\sin{\frac{\pi n x^2}{2}}\,dx$,$\displaystyle \ n>0$.

In analogous manner, one can deduce from \eqref{I2} and \eqref{int3} that \begin{equation}\label{factorization2}
    \int\limits_0^\infty \int\limits_0^\infty\frac{xye^{-x^2-y^2}\sin 2xy}{\cosh \alpha x\cosh({2\pi y}/{\alpha})}\,dxdy=\frac{\alpha}{2\sqrt{\pi}}\left(\int\limits_0^\infty\frac{\sinh \frac{\alpha x}{2}}{\cosh \alpha x}\,xe^{-x^2}dx \right)^2.
\end{equation}
${\it{Corollary}}~2.$
$$
\int\limits_0^\infty \int\limits_0^\infty\frac{\cos\frac{\pi}{2}
\left(nx^2-\frac{y^2}{n}\right)\sin \pi xy}{\cosh \pi x\cosh \pi y}\,xy\, dxdy=\frac{{\sqrt{n^3}}}{2}\left(I_4^2-I_3^2+2I_3I_4\right),
$$
$$
\int\limits_0^\infty \int\limits_0^\infty\frac{\sin\frac{\pi}{2}
\left(nx^2-\frac{y^2}{n}\right)\cos \pi xy}{\cosh \pi x\cosh \pi y}\,xy\, dxdy=\frac{{\sqrt{n^3}}}{2}\left(I_3^2-I_4^2+23_1I_4\right),
$$
where $\ \displaystyle I_3=\int\limits_0^\infty\frac{x\sinh \frac{\pi x}{2}}{\cosh \pi x}\cos{\frac{\pi n x^2}{2}}\,dx$,$\ \displaystyle I_4=\int\limits_0^\infty\frac{x\sinh \frac{\pi x}{2}}{\cosh \pi x}\sin{\frac{\pi n x^2}{2}}\,dx$,$\ \displaystyle n>0$.

Ramanujan showed that integrals $I_1-I_4$ have closed form expressions when $n\in\mathbb{Q}$ [\onlinecite{berndt}]. So the corresponding two-dimensional integrals also have closed form expressions.

$\it{Examples}.$
$$
\int\limits_0^\infty \int\limits_0^\infty\frac{\cos\frac{\pi}{2}\!
\left(3x^2-\frac{y^2}{3}\right)\cos \pi xy}{\cosh \pi x\cosh \pi y}\,dxdy=\frac{\sqrt{3}-1}{2\sqrt{6}},
$$
$$
\int\limits_0^\infty \int\limits_0^\infty\frac{\sin\frac{\pi}{2}\!
\left(3x^2-\frac{y^2}{3}\right)\cos \pi xy}{\cosh \pi x\cosh \pi y}\,dxdy=\frac{2-\sqrt{3}}{4\sqrt{2}},
$$
$$
\int\limits_0^\infty \int\limits_0^\infty\frac{\cos\frac{\pi}{2}\!
\left(x^2-{y^2}\right)\sin \pi xy}{\cosh \pi x\cosh \pi y}\,xy\, dxdy=\frac{1}{{8\sqrt{2}\pi^2}}.
$$

It is possible to calculate even more general integrals. In analogy with Ramanujan's integral analogs of theta functions [\onlinecite{berndt}] define 
\begin{equation}\label{phi1}
    \Phi_{\alpha,\beta}\left(\theta,\phi\right)=\int\limits_0^\infty \int\limits_0^\infty\frac{\cos \pi xy\cos\pi\theta x\cos\pi\phi y}{\cosh \pi x\cosh \pi y}\,e^{-\pi
(\alpha x^2+\beta y^2)/2}\,dxdy.
\end{equation}

One can apply the method developed by Ramanujan [\onlinecite{berndt}]  to the function $\Phi_{\alpha,\beta}\left(\theta,\phi\right)$.
From the definition of $\Phi_{\alpha,\beta}\left(\theta,\phi\right)$, it follows that
\begin{align}\label{2dmordellintegral2}
    \Phi_{\alpha,\beta}\left(\theta+i,\phi\right)+\Phi_{\alpha,\beta}\left(\theta-i,\phi\right)=
e^{-{\pi\theta^2}/(2\alpha)}\sqrt{\frac{2}{\alpha}}\int\limits_0^\infty\frac{\cos \pi\phi y \cosh \frac{\pi\theta y}{\alpha}}{\cosh \pi y }\,e^{{-{\pi}\left(\beta+1/{\alpha}\right)y^2/2}}\,dy.
\end{align}

Another transformation is derived similarly to (\ref{landen2})
\begin{align}\label{2dmordellintegral1}
\nonumber &\sqrt{\alpha\beta}\,e^{\pi
\theta^2/(2\alpha)+\pi\phi^2/(2\beta)}\,
\Phi_{\alpha,\beta}\left(\theta,\phi\right)+\Phi_{1/\alpha,1/\beta}({i\theta}/{\alpha},{i\phi}/{\beta})\\
&=\sqrt{2}\int\limits_0^\infty\frac{\cosh \frac{\pi x}{2}\cosh \frac{\pi\theta x}{\alpha} }{\cosh \pi x}\,e^{-\pi x^2/(2\alpha)}\,dx\cdot 
\int\limits_0^\infty\frac{\cosh \frac{\pi y}{2}\cosh \frac{\pi\phi y}{\beta} }{\cosh \pi y}\,e^{-\pi y^2/(2\beta)}\,dy.
\end{align}
Equation \eqref{2dmordellintegral1} generalizes \eqref{factorization1}: when $\alpha\beta=1$ and $i\theta/\alpha=\phi$ the double integrals in (\ref{2dmordellintegral1}) are equal and hence factorize.

Now we combine \eqref{2dmordellintegral1} and \eqref{2dmordellintegral2} to get
\begin{align}\label{2dmordellintegral3}
    &\sqrt{\frac{\beta}{2}}\left(e^{\pi\theta}\,\Phi_{\alpha,\beta}\left(\theta+\alpha,\phi\right)+e^{-\pi\theta}\,\Phi_{\alpha,\beta}\left(\theta-\alpha,\phi\right)\right)e^{\pi\phi^2/(2\beta)+\pi\alpha/2}\nonumber\\
&=-\int\limits_0^\infty\frac{\cos \pi\theta y \cosh \frac{\pi\phi y}{\beta}}{\cosh \pi y }\,e^{{-{\pi}\left({\alpha}+1/\beta\right)y^2/2}}dy+\sqrt{2}\,e^{{\pi\alpha}/{8}}\cosh\frac{\pi\theta}{2}\cdot\int\limits_0^\infty\frac{\cosh \frac{\pi y}{2}\cosh \frac{\pi\phi y}{\beta} }{\cosh \pi y}\,e^{-\pi y^2/(2\beta)}\,dy.
\end{align}
Thus when $\alpha/i\in\mathbb{Q}$ is a rational number, formulas (\ref{2dmordellintegral1}-\ref{2dmordellintegral3}) reduce the problem to the calculation of one-dimensional Mordell integrals. This shows that when $\alpha/i\in\mathbb{Q}$ and $\beta/i\in\mathbb{Q}$ are both rational, $\Phi_{\alpha,\beta}\left(\theta,\phi\right)$ can be calculated in closed form. Similar formulas exist for
$$
\Psi_{\alpha,\beta}\left(\theta,\phi\right)=\int\limits_0^\infty \int\limits_0^\infty\frac{\sin \pi xy\sin\pi\theta x\sin\pi\phi y}{\cosh \pi x\cosh \pi y}\,e^{-\pi
(\alpha x^2+\beta y^2)/2}\,dxdy.
$$

\section{Reduction of certain family of double Mordell integrals to combination of one-dimensional Mordell integrals}
Consider the following generalization of (\ref{phi1})
$$
\Phi_{\alpha,\beta}^{(\gamma)}=\int\limits_0^\infty \int\limits_0^\infty\frac{\cos \pi \gamma xy}{\cosh \pi x\cosh \pi y}\,e^{-\pi
(\alpha x^2+\beta y^2)/2}\,dxdy.
$$
We want to find a combination of parameters $\gamma$, $\alpha$, $\beta$ such that this integral reduces to a sum of products of one-dimensional Mordell integrals which we define according to Ramanujan [\onlinecite{berndt}] as
\begin{equation}\label{mordell}
    \phi_\alpha(\theta)=\int\limits_0^\infty\frac{\cos\pi \theta x}{\cosh \pi x}\,e^{-\pi \alpha x^2}dx.
\end{equation}
(\ref{mordell}) satisfies the transformation formula ([\onlinecite{berndt}], Entry 14.3.1)
\begin{equation*}
    \phi_\alpha(\theta)=\frac{1}{\sqrt{\alpha}}\,e^{-\pi\theta^2/(4\alpha)}\phi_{1/{\alpha}}(i\theta/\alpha).
\end{equation*}
First, we apply a series of transformations to $\Phi_{\alpha,\beta}^{(\gamma)}$:
\begin{align*}
\nonumber\Phi_{\alpha,\beta}^{(\gamma)}&=\int\limits_0^\infty\frac{e^{-\pi\beta y^2/2}}{\cosh \pi y}\,\phi_{\alpha/2}(\gamma y)\,dy\\
&=\sqrt{\frac{2}{\alpha}}\int\limits_0^\infty\frac{e^{{-{\pi}\left(\beta+\gamma^2/{\alpha}\right)y^2/2}}}{\cosh \pi y}\,\phi_{2/\alpha}(2i\gamma y/\alpha)\,dy\\
&=\sqrt{\frac{2}{\alpha}}\int\limits_{0}^\infty\int\limits_{0}^\infty\frac{e^{-2\pi x^2/\alpha-{\pi}\left(\beta+\gamma^2/{\alpha}\right)y^2/2}}{\cosh \pi x\,\cosh \pi y}\,\cosh\frac{2\pi \gamma xy}{\alpha}\,dxdy.
\end{align*}
It is convenient to extend integration over the whole plane:
\begin{align*}
\nonumber\Phi_{\alpha,\beta}^{(\gamma)}&=\sqrt{\frac{1}{8\alpha}}\int\limits_{-\infty}^\infty\int\limits_{-\infty}^\infty\frac{e^{-2\pi x^2/\alpha-2\pi \gamma xy/\alpha-{\pi}\left(\beta+\gamma^2/{\alpha}\right)y^2/2}}{\cosh \pi x\,\cosh \pi y}\,dxdy\\
&=\sqrt{\frac{1}{8\alpha}}\int\limits_{-\infty}^\infty\int\limits_{-\infty}^\infty\frac{e^{-2\pi (x+\gamma y/2)^2/\alpha-\pi\beta y^2/2}}{\cosh \pi x\,\cosh \pi y}\,dxdy\\
&=\sqrt{\frac{1}{8\alpha}}\int\limits_{-\infty}^\infty\int\limits_{-\infty}^\infty\frac{e^{-2\pi x^2/\alpha-\pi\beta y^2/2}}{\cosh \pi (x-\gamma y/2)\,\cosh \pi y}\,dxdy\\
&=\sqrt{\frac{1}{32\alpha}}\int\limits_{-\infty}^\infty\int\limits_{-\infty}^\infty\frac{e^{-2\pi x^2/\alpha-\pi\beta y^2/2}}{\cosh \pi y}\left(\frac{1}{\cosh \pi (x-\gamma y/2)}+\frac{1}{\cosh \pi (x+\gamma y/2)}\right)dxdy\\
&=\sqrt{\frac{1}{8\alpha}}\int\limits_{-\infty}^\infty\int\limits_{-\infty}^\infty\frac{e^{-2\pi x^2/\alpha-\pi\beta y^2/2}\cosh\pi x}{\cosh \pi (x-\gamma y/2)\cosh \pi (x+\gamma y/2)}\frac{\cosh\frac{\pi \gamma y}{2}}{\cosh\pi y}\,dxdy.
\end{align*}
After the change of variables $\xi=x-\gamma y/2,~\eta=x+\gamma y/2$, considerable simplification occurs when $\alpha \beta=\gamma^2$:
\begin{align*}
\nonumber\Phi_{\alpha,\gamma^2/\alpha}^{(\gamma)}&=\frac1n\sqrt{\frac{1}{8\alpha}}\int\limits_{-\infty}^\infty\int\limits_{-\infty}^\infty\frac{e^{-\pi (\xi^2+\eta^2)/\alpha}\cosh\frac{\pi(\xi+\eta)}{2}}{\cosh \pi\xi\cosh \pi\eta}\frac{\cosh\frac{\pi(\xi-\eta)}{2}}{\cosh\frac{\pi(\xi-\eta)}{\gamma}}\,dxdy.
\end{align*}
To complete the process of reduction we set $\gamma=4n+2$, $n\in\mathbb{N}_0$:
\begin{align*}
\nonumber\Phi_{\alpha,(4n+2)^2/\alpha}^{(4n+2)}&=\frac{1}{4n+2}\sqrt{\frac{1}{8\alpha}}\int\limits_{-\infty}^\infty\int\limits_{-\infty}^\infty\frac{e^{-\pi (\xi^2+\eta^2)/\alpha}\cosh\frac{\pi(\xi+\eta)}{2}}{\cosh \pi\xi\cosh \pi\eta}\sum_{k=-n}^n(-1)^{k}e^{\pi k(\xi-\eta)/(2n+1)}\,dxdy\\
&=\frac{1}{2n+1}\sqrt{\frac{1}{2\alpha}}\int\limits_{0}^\infty\int\limits_{0}^\infty\frac{e^{-\pi (\xi^2+\eta^2)/\alpha}}{\cosh \pi\xi\cosh \pi\eta}\sum_{k=-n}^n(-1)^{k}\cosh\tfrac{\pi(2n+2k+1)\xi}{4n+2}\cosh\tfrac{\pi(2n-2k+1)\xi}{4n+2}\,dxdy.
\end{align*}
Thus, we have proved the first reduction formula
\begin{equation}\label{reduction1}
    \Phi_{1/\alpha, (4n+2)^2\alpha}^{(4n+2)}=\frac{1}{2n+1}\sqrt{\frac{\alpha}{2}}\left(\left\{\phi_{\alpha}\left(\tfrac{i}{2}\right)\right\}^2+2\sum_{k=1}^n(-1)^{k}\phi_{\alpha}\Big(\tfrac{2n+2k+1}{4n+2}i\Big)\phi_{\alpha}\Big(\tfrac{2n-2k+1}{4n+2}i\Big)\right),\quad n\in\mathbb{N}_0.
\end{equation}

There is a transformation formula between two functions $\Phi_{\alpha,\beta}^{(\gamma)}$ that we will now derive. First, iterating the gaussian integral  $4.133.2$ from [\onlinecite{GR}]
$$
\int\limits_0^\infty e^{-{x^2}/{(4c)}}\cos ax\cosh bx\, dx=\sqrt{\pi c}\, e^{c(b^2-a^2)}\cos (2abc),
$$
where $\text{Re}\,c>0$, one can show that

$$
\int\limits_0^\infty \int\limits_0^\infty e^{-(x^2+y^2)/2}\cos(qxy)\cos ax\cos bx\, dxdy=\frac{\pi}{2\sqrt{1+q^2}}\, \exp\left\{{-\frac{a^2+b^2}{2(1+q^2)}}\right\}\cos \frac{qab}{1+q^2}.
$$
Then multiplying this integral by $1/\bigl(\cosh(\sqrt{\frac{\pi}{2}}\frac{ a}{\alpha})\cosh(\sqrt{\frac{\pi}{2}}\frac{ b}{\beta})\bigr)$ and integrating wrt $a$ and $b$ we come to
\begin{align*}
    \alpha\beta\int\limits_0^\infty \int\limits_0^\infty \frac{e^{-(x^2+y^2)/2}\,\cos(qxy)}{\cosh(\sqrt{\frac{\pi}{2}}{\alpha}x)\cosh(\sqrt{\frac{\pi}{2}}{\beta}y)}\,dxdy=\frac{1}{\sqrt{1+q^2}}\int\limits_0^\infty \int\limits_0^\infty \exp\left\{{-\frac{x^2+y^2}{2(1+q^2)}}\right\}\frac{\cos \frac{qxy}{1+q^2}\,dxdy}{\cosh(\sqrt{\frac{\pi}{2}}\frac{ x}{\alpha})\cosh(\sqrt{\frac{\pi}{2}}\frac{y}{\beta})}.
\end{align*}
This implies the following general three-parameter transformation for $\Phi_{\alpha,\beta}^{(\gamma)}$:
\begin{equation*}
    \Phi_{{2}/{\alpha},{2}/{\beta}}^{(2q/\sqrt{\alpha\beta})}=\sqrt{\frac{\alpha\beta}{1+q^2}}\,\Phi_{{2\alpha}/{(1+q^2)},{2\beta}/{(1+q^2)}}^{(2q\sqrt{\alpha\beta}/(1+q^2))}.
\end{equation*}
Combining with (\ref{reduction1}) we find another family of double Mordell integrals that reduce to a combination of one-dimensional integrals
\begin{equation}\label{reduction2}
    \sqrt{(4n+2)\alpha}\cdot \Phi_{1/\alpha,\alpha/(2n+1)^2}^{(1/(2n+1))}=\left\{\phi_{1/(2\alpha)}\left(\tfrac{i}{2}\right)\right\}^2+2\sum_{k=1}^n(-1)^{k}\phi_{1/(2\alpha)}\Bigl(\tfrac{2n+2k+1}{4n+2}i\Bigr)\phi_{1/(2\alpha)}\Bigl(\tfrac{2n-2k+1}{4n+2}i\Bigr),\quad n\in\mathbb{N}_0.
\end{equation}
This is second reduction formula. (\ref{reduction1}) and (\ref{reduction2}) are main formulas of this section. Note that both double Mordell integrals in (\ref{reduction1}) and (\ref{reduction2}) are of the type $\Phi_{\alpha,\beta}^{(\sqrt{\alpha\beta})}$.

Two examples of the reduction formula (\ref{reduction2}) are shown below:

i) $n=0$: In this case we recover (\ref{landen2}).

ii) $n=1$: 
\begin{align*}
\sqrt{\frac{3}{\alpha}}\int\limits_0^\infty\int\limits_0^\infty &\frac{\cos(\pi xy/3)}{\cosh \pi x\cosh \pi y}\,e^{-\pi
\alpha x^2-\pi y^2/(36\alpha)}\,dxdy\\&=\bigg(\int\limits_0^\infty\frac{\cosh \frac{\pi  x}{2} }{\cosh \pi  x}\,e^{-\pi  \alpha x^2}dx\bigg)^2-2\int\limits_0^\infty\frac{\cosh \frac{\pi  x}{6} }{\cosh \pi  x}\,e^{-\pi  \alpha x^2}dx\int\limits_0^\infty\frac{\cosh \frac{5\pi  x}{6} }{\cosh \pi  x}\,e^{-\pi  \alpha x^2}dx.
\end{align*}

There is a curious consequence of the formulas above for one-dimensional Mordell integrals. For $n=0$ it is just equation (\ref{int2}). For $n=1$ it reads as follows:
$$
\sqrt{\alpha}\left\{\phi_{\alpha}\left(\tfrac{i}{2}\right)\right\}^2-2\sqrt{\alpha}\,\phi_{\alpha}\big(\tfrac{i}{6}\big)\phi_{\alpha}\big(\tfrac{5i}{6}\big)=
\sqrt{\beta}\left\{\phi_{\beta}\left(\tfrac{i}{2}\right)\right\}^2-2\sqrt{\beta}\,\phi_{\beta}\big(\tfrac{i}{6}\big)\phi_{\beta}\big(\tfrac{5i}{6}\big),\qquad \alpha\beta=1/36.
$$
This is a quadratic relation connecting $6$ different Mordell integrals. Linear relations between one-dimensional Mordell integrals have been studied before (e.g., [\onlinecite{watson2}]) and two-dimensional Mordell integrals have been investigated recently in connection with vector-valued higher depth quantum modular forms [\onlinecite{bringmann}]. However, it seems the fact that there are non-trivial reductions of certain two-dimensional Mordell integrals to one-dimensional Mordell integrals, or the fact that there are non-trivial quadratic relations between one-dimensional Mordell integrals have not been recognized in the existing literature.

\section{Absolute value of the Mordell integral}
In this section we study integrals of the type
$$
\int\limits_0^\infty\frac{e^{i\alpha x^2}}{\cosh\pi x}\,dx,
$$
where $\alpha\in\mathbb{R}$. The square of the absolute value of this integral can be transformed in the following way:
\begin{align*}
4\Bigg|\int\limits_0^\infty\frac{e^{i\alpha x^2}}{\cosh\pi x}\,dx\Bigg|^2&=\int\limits_{-\infty}^\infty\frac{e^{i\alpha x^2}}{\cosh\pi x}\,dx\int\limits_{-\infty}^\infty\frac{e^{-i\alpha (x+y)^2}}{\cosh\pi(x+y)}\,dy\\
&=\int\limits_{-\infty}^\infty \int\limits_{-\infty}^\infty\frac{e^{-i\alpha y^2-2i\alpha xy}}{\cosh\pi x\cosh\pi(x+y)}\,dxdy\\
&=\int\limits_{-\infty}^\infty \int\limits_{-\infty}^\infty\frac{e^{-2i\alpha xy}}{\cosh\pi (x-y/2)\cosh\pi(x+y/2)}\,dxdy\\
&=2\int\limits_{-\infty}^\infty dy \int\limits_{-\infty}^\infty\frac{e^{-2i\alpha xy}}{\cosh2\pi x+\cosh\pi y}\,dx\\
&=2\int\limits_{-\infty}^\infty\frac{\sin\alpha y^2}{\sinh\pi y\sinh\alpha y}\,dy.
\end{align*}
This can be written as
\begin{equation}
    \int\limits_0^\infty\frac{\sin\alpha x^2}{\sinh\pi x\sinh\alpha x}\,dx=\Bigg(\int\limits_0^\infty\frac{\cos{\alpha x^2}}{\cosh\pi x}\,dx\Bigg)^2+\Bigg(\int\limits_0^\infty\frac{\sin{\alpha x^2}}{\cosh\pi x}\,dx\Bigg)^2, \qquad \alpha\in\mathbb{R}.
\end{equation}

Analogous considerations lead to other formulas of similar kind
\begin{equation}\label{one_half}
    \int\limits_0^\infty\frac{\sin 2\alpha x^2}{\sinh\pi x\sinh \alpha x}\,dx=\int\limits_0^\infty\frac{\cos 2\alpha x^2}{\cosh\pi x\cosh \alpha x}\,dx=\Bigg|\int\limits_0^\infty\frac{\cosh\frac{\pi x}{2}}{\cosh\pi x}\,e^{{i \alpha x^2}/{2}}\,dx\Bigg|^2,\qquad \alpha\in\mathbb{R},
\end{equation}
\begin{equation}
    \pi\int\limits_0^\infty\frac{\sin\frac{3\alpha x^2}{4\pi}\coth\frac{x}{2}\coth\frac{\alpha x}{2}-\frac{1}{\sqrt{3}}\cos \frac{3\alpha x^2}{4\pi}}{(1+2\cosh x)(1+2\cosh\alpha x)}\,dx=\Bigg|\int\limits_0^\infty\frac{e^{{ 3i\alpha x^2}/(4\pi)}}{1+2\cosh x}\,dx\Bigg|^2,\qquad \alpha\in\mathbb{R},
\end{equation}
and to the following curious closed form
\begin{equation}\label{closed_form}
    \int\limits_0^\infty\tanh {\pi x}\tanh{\alpha x}\cos{2\alpha x^2}\,dx=0,\qquad \alpha\in\mathbb{R}.
\end{equation}

Here we give an explanation for the first equality in \eqref{one_half} and for \eqref{closed_form}. For the first, starting from \eqref{f1} we put $b=\alpha a$
\begin{equation*}
    {\frac{2}{\pi}}\int\limits_0^\infty \int\limits_0^\infty \frac{\cos xy}{\cosh\sqrt{\frac{\pi}{2}} x\cdot \cosh\sqrt{\frac{\pi}{2}} y}\cos ax\cos \alpha ay\,dxdy=\frac{\sin \alpha a^2}{\sinh\sqrt{\frac{\pi}{2}} a \cdot\sinh\sqrt{\frac{\pi}{2}}\alpha a}
\end{equation*}
and integrate with respect to $a$ from $0$ to $\infty$ 
to obtain
$$
\int\limits_0^\infty\frac{\cos 2\alpha x^2}{\cosh\pi x\cosh \alpha x}\,dx=\int\limits_0^\infty\frac{\sin 2\alpha x^2}{\sinh\pi x\sinh \alpha x}\,dx.
$$
For the second, starting from \eqref{cos2} and its sine analog
\begin{equation}\label{sine_analog2}
    \frac{2}{\pi}\int\limits_0^\infty \int\limits_0^\infty\frac{\cos xy}{\sinh\sqrt{\pi} x\, \sinh\sqrt{\pi}y}\sin ax\sin by\,dxdy=\frac12\tanh\frac{\sqrt{\pi} a}{2}\tanh\frac{\sqrt{\pi} b}{2}-\frac{1-\cos ab}{\sinh\sqrt{\pi} a\, \sinh\sqrt{\pi}b},
\end{equation}
(by the way, \eqref{sine_analog2} implies the self reciprocal function $\frac{1-\cos xy}{\sinh\sqrt{\pi} x\, \sinh\sqrt{\pi}y}$)
we put $b=\alpha a$ in both, take the sum of \eqref{cos2} multiplied by $e^{{i\alpha} a^2/2}$ and \eqref{sine_analog2} multiplied by $ie^{{i\alpha} a^2/2}$, integrate from $0$ to $\infty$ 
using formulas
$$
\int\limits_0^\infty \cos ax\cos\alpha ay\,e^{{i\alpha} a^2/2}\,da=\sqrt{\frac{\pi i}{2\alpha}}\,e^{-{i}(x^2+\alpha^2 y^2)/(2\alpha)}\cos xy,
$$
$$
\int\limits_0^\infty \sin ax\sin\alpha ay\,e^{{i\alpha} a^2/2}\,da=i\sqrt{\frac{\pi i}{2\alpha}}\,e^{-{i}(x^2+\alpha^2 y^2)/(2\alpha)}\sin xy,
$$
to obtain
$$
0=\frac{i}2\int\limits_0^\infty\tanh\frac{\sqrt{\pi} a}{2}\tanh\frac{\sqrt{\pi} \alpha a}{2}\,e^{{i\alpha} a^2/2}\,da+\int\limits_0^\infty \frac{-i(1-\cos \alpha a^2)+\sin\alpha a^2}{\sinh\sqrt{\pi} a \sinh\sqrt{\pi}\alpha a}\,e^{{i\alpha} a^2/2}\,da.
$$
From this, it is straightforward to deduce \eqref{closed_form} and as a byproduct
$$
\int\limits_0^\infty\frac{2\sin\frac{\alpha x^2}{2}}{\sinh \pi x\sinh\alpha x}\,dx=\int\limits_0^\infty\tanh {\pi x}\tanh{\alpha x}\sin{2\alpha x^2}\,dx, \qquad \alpha\in\mathbb{R}.
$$

Note the equivalent formulation of (\ref{closed_form}):
\begin{equation}\label{closed_form2}
    \int\limits_0^\infty\frac{\cosh {(\pi-\alpha) x}}{\cosh\pi x\cosh{\alpha x}}\cos{2\alpha x^2}\,dx=\frac{1}{4}\sqrt{\frac{\pi}{\alpha}},\qquad \alpha>0.
\end{equation}
Formulas \eqref{closed_form} and \eqref{closed_form2} are reminiscent of the integral of Ramanujan 
\begin{equation}\label{closed_form3}
    \int\limits_0^\infty\frac{\cosh \alpha x}{\cosh \pi x}\cos{\alpha x^2}\,dx=\frac{1}{2}\cos\frac{\alpha}{4}, \qquad \alpha\in\mathbb{R},
\end{equation}
([\onlinecite{ramanujan2}], see also generalizations in [\onlinecite{glasser1},\onlinecite{glasser2}]). \eqref{closed_form}, \eqref{closed_form2} and \eqref{closed_form3} contain trigonometric function of the argument ${\alpha x^2}$ and hyperbolic functions of the arguments $\pi x$ and $\alpha x$. However the crucial difference between them is that the integrand in \eqref{closed_form} has poles not only at the zeroes of $\cosh \pi x$, but also at the zeroes of $\cosh \alpha x$. Integrals of this sort are related to integrals for the product of two hyperbolic self-reciprocal functions studied by Ramanujan ([\onlinecite{ramanujan}], formula $(10)$). To show this we put $b=\alpha a$ in \eqref{I1} and \eqref{I2} and integrate with respect to $a$. The result is
\begin{equation}
   \sqrt{2}\int\limits_0^\infty \frac{\cos \alpha x^2}{\cosh \pi  x \cosh \alpha x}\,dx=\int\limits_0^\infty  \frac{\cosh \frac{\pi  x}{2}}{\cosh \pi  x }\cdot\frac{ \cosh \frac{\alpha x}{2}}{\cosh \alpha x}\, dx, \qquad \alpha\in\mathbb{R},
\end{equation}
\begin{equation}
   \sqrt{2}\int\limits_0^\infty \frac{\sin \alpha x^2}{\cosh \pi  x \cosh \alpha x}\,dx=\int\limits_0^\infty \frac{\sinh \frac{\pi  x}{2}}{\cosh \pi  x }\cdot\frac{ \sinh \frac{\alpha x}{2}}{\cosh \alpha x}\, dx, \qquad \alpha\in\mathbb{R}.
\end{equation}

\section{Connection to lattice sums}

Multiplying \eqref{I1} and \eqref{I2} by $\displaystyle\frac{1}{\sqrt{ab}}$ and integrating with respect to $a$ and $b$ leads to
\begin{equation}\label{lattice1}
    \sqrt{2}\int\limits_0^\infty \int\limits_0^\infty\frac{\cos\frac{x^2y^2}{\pi}\,dxdy}{\cosh x^2\cosh y^2}=\left(\int\limits_0^\infty\frac{\cosh\frac{x^2}{2}}{\cosh x^2}\,dx\right)^2,
\end{equation}
\begin{equation}\label{lattice2}
    \sqrt{2}\int\limits_0^\infty \int\limits_0^\infty\frac{\sin\frac{x^2y^2}{\pi}\,dxdy}{\cosh x^2\cosh y^2}=\left(\int\limits_0^\infty\frac{\sinh\frac{x^2}{2}}{\cosh x^2}\,dx\right)^2.
\end{equation}
The RHS of \eqref{lattice1} and \eqref{lattice2} contain integral representation of certain Dirichlet L-series, while the LHS are 2D-lattice sums of Bessel and Neumann functions, as shown below on a formal level. Evaluation of double sums of Bessel functions in terms of Dirichlet L-series is well known [\onlinecite{lattice}].

Consider the double integral on the LHS of \eqref{lattice1}. First, the functions $\text{sech} ~x^2$ are expanded into the powers of $e^{-x^2}$. This results in a double sum of double integrals
$$
\int\limits_0^\infty \int\limits_0^\infty e^{-(2m+1)x^2-(2n+1)y^2}\cos\frac{x^2y^2}{\pi}~dxdy,
$$
where $m$ and $n$ are non-negative integers. The integral over $y$ is easily calculated
$$
\int\limits_0^\infty e^{-(2n+1)y^2}\cos\frac{x^2y^2}{\pi}\,dy=\frac{\pi}{2}\left(\frac{1}{\sqrt{\pi(2m+1)+ix^2}}+\frac{1}{\sqrt{\pi(2m+1)-ix^2}}\right).
$$
To calculate the integral over $x$ we need formula $3.364.3$ from [\onlinecite{GR}]
$$
\int\limits_0^\infty \frac{e^{-(2n+1)x^2}}{\sqrt{\pi(2m+1)\pm ix^2}}\,dx=\frac12 (-1)^{m+n}e^{\mp \frac{3\pi i}{4}}K_0\left(\mp \frac{\pi i}{2}(2m+1)(2n+1))\right).
$$
Note that $K_0(ix)=-\frac{\pi}{2}(Y_0(x)+iJ_0(x))$,$~x\in\mathbb{R}$. As a result the double integral in \eqref{lattice1} reduces to a combination of double sums
$$
\sum_{m,n=0}^\infty Z_0\left(\frac{\pi}{2}(2m+1)(2n+1)\right),
$$
where $Z_0$ is either Bessel $J_0$ or Neumann $Y_0$ function.

{\it{Acknowledgements.}} The author of this paper wish to thank Dr. Lawrence Glasser for valuable correspondence and comments.


\begin{thebibliography}{9}
\bibitem{berndt} B. Berndt and G.E. Andrews, {\it{Ramanujan's lost notebook, part IV}}, Springer New York (2005).

\bibitem{lattice}  J.M. Borwein, M.L. Glasser, R.C. McPhedran, J.G. Wan, I.J. Zucker, {\it{Lattice sums then and now}}, Cambridge University Press  (2013).
\bibitem{bringmann} K. Bringmann, J. Kaszian, A. Milas, {\it{Vector-valued higher depth quantum modular forms and higher Mordell integrals}}, Journal of Mathematical Analysis and Applications, Volume 480, Issue 2, p. 123397 (2019).
\bibitem{glasser1} M. L. Glasser, {\it{Generalization of a definite integral of Ramanujan}}, J. Indian Math. Soc., 37, 351 (1973).
\bibitem{glasser2} M. L. Glasser, {\it{A Remarkable Definite Integral}},  {\href{http://arxiv.org/abs/1308.6361}{arXiv:1308.6361v2}} (2013).
\bibitem{GR} I.S. Gradshteyn, and I.M. Ryzhik, {\it{Table of Integrals, Series, and Products}}, 6th ed.,  Academic Press, Boston (2000).
\bibitem{ramanujan2} S. Ramanujan, {\it{Some definite integrals}}, J. Indian Math. Soc., XI, 81-87 (1919).
\bibitem{hardy} G.H. Hardy, {\it{Note on the function}$~\int_x^\infty e^{\frac12(x^2-t^2)}dt$}, Quart. J. Pure Appl. Math., 35, 203 (1903).
\bibitem{cais} B. Cais, {\it{On the transformation of infinite series}} (unpublished), (1999).
\bibitem{ramanujan} S. Ramanujan, {\it{Some definite integrals}}, Mess. Math., XLIV, 10 - 18 (1915).
\bibitem{titchmarsh} E.C. Titchmarsh, {\it{Introduction to the Theory of Fourier Integrals}},  2nd.ed., Oxford University Press (1948).
\bibitem{watson2} G. N. Watson, {\it{Generating functions of class-numbers}}, Compositio Mathematica, v. 1, p. 39-68, (1935).
\bibitem{ww} E.T. Whittaker and G.N. Watson, {\it{A Course of Modern Analysis}}, Cambridge university press (1996).
\bibitem{blog} \url{www.someformulas.blogspot.com}




\end{thebibliography}
\end{document}